\documentclass[12pt]{article}
\usepackage{graphicx}
\usepackage{amsmath}
\usepackage{amsfonts}
\usepackage{amssymb}
\newtheorem{theorem}{Theorem}

\newtheorem{corollary}[theorem]{Corollary}

\newtheorem{lemma}[theorem]{Lemma}

\newtheorem{proposition}[theorem]{Proposition}

\newenvironment{proof}[1][Proof]{\textbf{#1.} }{\ \rule{0.5em}{0.5em}}

\def\text{\hbox} 

\def\a{\alpha}
\def\b{\beta}
\def\g{\gamma}

\def\m{\mu}

\def\p{\pi}

\def\r{\rho}

\def\s{\sigma}

\def\t{\tau}

\def\f{\phi}

\def\ps{\psi}

\def\o{\omega}

\def\Si{\Sigma}

\def\F{\Phi}
\def\PS{\Psi}

\def\S{{\bf S}}
\def\E{{\bf E}}
\def\P{{\cal P}}
\def\R{{\cal R}}
\def\A{{\cal A}}
\def\B{{\cal B}}
\def\A'{{\cal A}'}
\def\B'{{\cal B}'}
\def\GG{{\cal G}}
\def\HH{{\cal H}}
\def\FF{{\cal F}}
\def\Z{{\mathbf Z}}
\def\wh{\widehat}

\def\wti{\widetilde}

\begin{document}

\title{On $p$-symmetric Heegaard splittings
\footnote{Work performed under the auspices of
G.N.S.A.G.A. of C.N.R. of Italy and supported by the University of
Bologna, funds for selected research topics.}}
\author{Michele Mulazzani}

\maketitle
\begin{abstract}
{We show that every $p$-fold strictly-cyclic branched covering of
a $b$-bridge link in $\S^3$ admits a $p$-symmetric Heegaard
splitting -- in the sense of Birman and Hilden -- of genus
$g=(b-1)(p-1)$. This gives a complete converse of one of the
results of the two authors. Moreover, we introduce the concept of
weakly $p$-symmetric Heegaard splittings and prove similar results
in this more general context.
}

\end{abstract}

%\newpage

\section{Introduction}

Let $M$ be a closed, orientable 3-manifold, $Y_g$ and $Y'_g$ be
handlebodies of genus $g\ge 0$ and $\f:\partial Y_g\to\partial
Y'_g$ be an attaching homeomorphism. A Heegaard splitting
$M=Y_g\cup_{\f}Y'_g$ is called in \cite {BH} {\it
$p$-symmetric\/}, with $p>1$, if there is a disjoint embedding of
$Y_g$ and $Y'_g$ into $\E^3$ such that $Y'_g=\t(Y_g)$, for a
translation $\t$ of $\E^3$, and an orientation-preserving
homeomorphism $\P:\E^3\to\E^3$ of period $p$, such that
$\P(Y_g)=Y_g$ and, if $\GG$ is the cyclic group of order $p$
generated by $\P$ and $\F:\partial Y_g\to\partial Y_g$ is the
orientation-preserving homeomorphism $\F=\t^{-1}_{\vert\partial
Y'_g}\f$, the following conditions are fulfilled:
\begin{itemize}
\item[i)] $Y_g/\GG$ is homeomorphic to a 3-ball;
\item[ii)] $\mbox{Fix}(\P_{\vert Y_g}^h)=\mbox{Fix}(\P_{\vert Y_g})$,
for each $1\le h\le p-1$;
\item[iii)] $\mbox{Fix}(\P_{\vert Y_g})/\GG$ is an unknotted set of
arcs\footnote{A set of mutually disjoint arcs $\{t_1,\ldots,t_n\}$
properly embedded in a handlebody $Y$ is {\it unknotted\/} if
there is a set of mutually disjoint discs $D=\{D_1,\ldots,D_n\}$
properly embedded in $Y$ such that $t_i\cap D_i=t_i\cap\partial
D_i=t_i$, $t_i\cap D_j=\emptyset$ and $\partial
D_i-t_i\subset\partial Y$ for $1\le i,j\le n$ and $i\neq j$.} in
the ball $Y_g/{\cal G}$;
\item[iv)] there exists an integer $p_0$ such that
$\F\P_{\vert\partial Y_g}\F^{-1}=(\P_{\vert\partial
Y_g})^{p_0}$.\footnote{By the positive solution of the Smith
Conjecture \cite{MB} it is easy to see that necessarily
$p_0\equiv\pm 1$ mod $p$.}
\end{itemize}

Observe that the map $\P'=\t\P\t^{-1}$ is obviously an
orientation-preserving homeomorphism of period $p$ of $\E^3$ with
the same properties as $\P$, regarding $Y'_g$, and the relation
$\f\P_{\vert\partial Y_g}\f^{-1}=(\P'_{\vert\partial Y'_g})^{p_0}$
easily holds.

The $p$-symmetric Heegaard genus $g_p(M)$ of a 3-manifold $M$ is
the smallest integer $g$ such that $M$ admits a $p$-symmetric
Heegaard splitting of genus $g$. The following results have been
established in \cite{BH}:

\begin{proposition} \label{Theorem 2BH} (\cite{BH}, Theorem 2)
Every closed, orientable 3-manifold of $p$-symmetric Heegaard
genus $g$ admits a representation as a $p$-fold cyclic covering of
$\S^3$, branched over a link $L$ of bridge number $$b(L)\le
1+\frac{g}{p-1}.$$
\end{proposition}

\begin{proposition} \label{Theorem 3BH} (\cite{BH}, Theorem 3) The
$p$-fold cyclic covering of $\S^3$ branched over a knot of braid
number $b$ is a closed, orientable 3-manifold $M$ of $p$-symmetric
Heegaard genus $$g_p(M)\le(b-1)(p-1).\eqno (*)$$
\end{proposition}

As a consequence, Birman and Hilden got the interesting result
that every closed, orientable 3-manifold of Heegaard genus $g\le
2$ is a 2-fold covering of $\S^3$ branched over a link of bridge
number $g+1$ and, conversely, the 2-fold covering of $\S^3$
branched over a link of bridge number $b\le 3$ is a closed,
orientable 3-manifold of Heegaard genus $b-1$ (compare also
\cite{Vi}).

Note that Proposition \ref{Theorem 3BH} is not a complete converse
of Proposition \ref{Theorem 2BH}, since it only concerns knots
and, moreover, the relation $(*)$ connects the genus with the
braid number, which is greater than or equal to (often greater)
the bridge number. In this paper we fill this gap, giving a
complete converse of Proposition \ref{Theorem 2BH}. Since the
coverings involved in Proposition \ref{Theorem 2BH} are
strictly-cyclic (see next chapter for details on strictly-cyclic
branched coverings of links), our statement will concern these
kind of coverings. More precisely, we shall prove in Theorem
\ref{Theorem 3} that a $p$-fold strictly-cyclic covering of
$\S^3$, branched over a link of bridge number $b$, is a closed,
orientable 3-manifold $M$ which admits a $p$-symmetric Heegaard
splitting of genus $g=(b-1)(p-1)$, and therefore has $p$-symmetric
Heegaard genus $g_p(M)\le (b-1)(p-1)$. The proof of this fact will
be achieved using the concept of special plat presentation of a
link.

In the last section we introduce a more intrinsic notion of
$p$-symmetric Heegaard splitting of a 3-manifold -- which will be
called weakly $p$-symmetric Heegaard splitting -- generalizing the
one given by Birman and Hilden. Proposition \ref{Theorem 2BH} will
translate to this new context and a complete converse of it will
be given, for all cases of branched cyclic coverings of links.

\section{Special plat presentations of links}

Let $\b=\{(p_k(t),t)\,\vert\, 1\le k\le
2n\,,\,t\in[0,1]\}\subset\E^2\times[0,1]$ be a geometric
$2n$-string braid of $\E^3$ \cite{Bi}, where
$p_1,\ldots,p_{2n}:[0,1]\to\E^2$ are continuous maps such that
$p_{k}(t)\neq p_{k'}(t)$, for every $k\neq k'$ and $t\in[0,1]$,
and such that
$\{p_1(0),\ldots,p_{2n}(0)\}=\{p_1(1),\ldots,p_{2n}(1)\}$. We set
$P_k=p_k(0)$, for each $k=1,\ldots,2n$, and
$A_i=(P_{2i-1},0),B_i=(P_{2i},0),A'_i=(P_{2i-1},1),B'_i=(P_{2i},1)$,
for each $i=1,\ldots,n$ (see Figure 1). Moreover, we set
$\FF=\{P_1,\ldots,P_{2n}\}$, $\FF_1=\{P_1,P_3\ldots,P_{2n-1}\}$
and $\FF_2=\{P_2,P_4,\ldots,P_{2n}\}$.

The braid $\b$ is realized through an ambient isotopy
${\wh\b}:\E^2\times[0,1]\to\E^2\times[0,1]$,
${\wh\b}(x,t)=(\b_t(x),t)$, where $\b_t$ is an homeomorphism of
$\E^2$ such that $\b_0=\mbox{Id}_{\E^2}$ and $\b_t(P_i)=p_i(t)$,
for every $t\in[0,1]$. Therefore, the braid $\b$ naturally defines
an orientation-preserving homeomorphism
${\wti\b}=\b_1:\E^2\to\E^2$, which fixes the set $\FF$. Note that
$\b$ uniquely defines ${\wti\b}$, up to isotopy of $\E^2$ mod
$\FF$.

Connecting the point $A_i$ with $B_i$ by a circular arc $\a_i$
(called {\it top arc\/}) and the point $A'_i$ with $B'_i$ by a
circular arc $\a'_i$ (called {\it bottom arc\/}), as in Figure 1,
for each $i=1,\ldots,n$, we obtain a $2n$-plat presentation of a
link $L$ in $\E^3$, or equivalently in $\S^3$. As is well known,
every link admits plat presentations and, moreover, a $2n$-plat
presentation corresponds to an $n$-bridge presentation of the
link. So, the bridge number $b(L)$ of a link $L$ is the smallest
positive integer $n$ such that $L$ admits a representation by a
$2n$-plat. For further details on braids, plat and bridge
presentations of links we refer to \cite{Bi}.

%\bigskip

\begin{figure}[bht]
 \begin{center}
 \includegraphics*[totalheight=6cm]{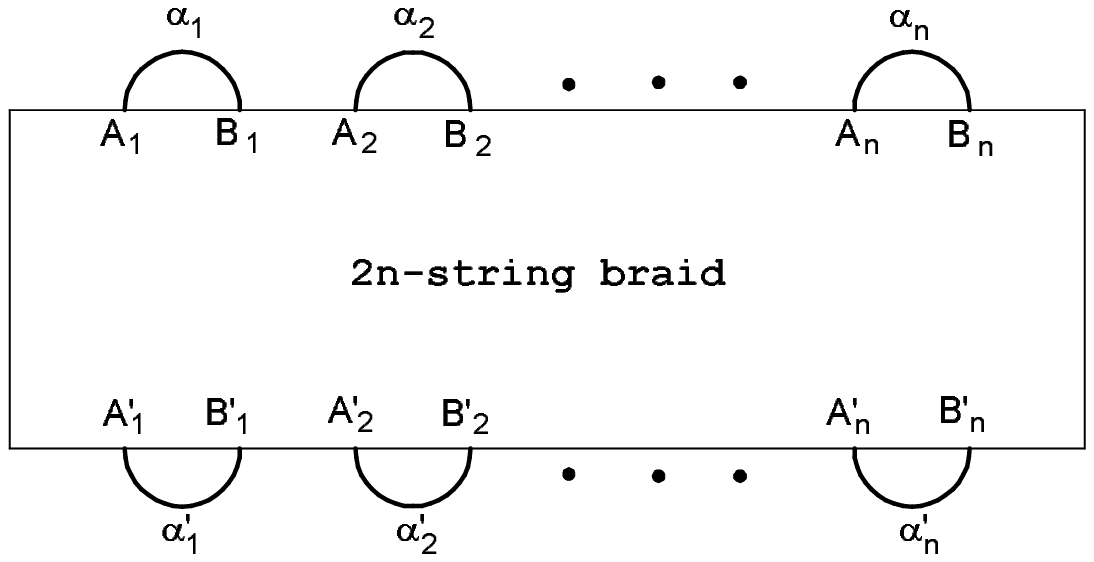}
 \end{center}
 \caption{A $2n$-plat presentation of a link.}

 \label{Fig. 1}

\end{figure}

\medskip

\noindent {\bf Remark 1} A $2n$-plat presentation of a link
$L\subset\E^3\subset\S^3=\E^3\cup\{\infty\}$ furnishes a
$(0,n)$-decomposition \cite{MS}
$(\S^3,L)=(D,A_n)\cup_{\f'}(D',A'_n)$ of the link, where $D$ and
$D'$ are the 3-balls $D=(\E^2\times]-\infty,0])\cup\{\infty\}$ and
$D'=(\E^2\times[1,+\infty[)\cup\{\infty\}$,
$A_n=\a_1\cup\cdots\cup\a_n$, $A'_n=\a'_1\cup\cdots\cup\a'_n$ and
$\f':\partial D\to\partial D'$ is defined by $\f'(\infty)=\infty$
and $\f'(x,0)=({\wti\b}(x),1)$, for each $x\in\E^2$.

\medskip

If a $2n$-plat presentation of a $\m$-component link
$L=\bigcup_{j=1}^{\m}L_j$ is given, each component $L_j$ of $L$
contains $n_j$ top arcs and $n_j$ bottom arcs. Obviously,
$\sum_{j=1}^{\m}n_j=n$. A $2n$-plat presentation of a link $L$
will be called {\it special \/} if:
\begin{itemize}
\item[(1)] the top arcs and the bottom arcs belonging to
$L_1$ are $\a_1,\ldots,\a_{n_1}$ and $\a'_1,\ldots,\a'_{n_1}$, the
top arcs and the bottom arcs belonging to $L_2$ are
$\a_{n_1+1},\ldots,\a_{n_1+n_2}$ and
$\a'_{n_1+1},\ldots,\a'_{n_1+n_2}$, $\ldots$ , the top arcs and
the bottom arcs belonging on $L_\m$ are
$\a_{n_1+\cdots+n_{\m-1}+1},\ldots,\a_{n_1+\cdots+n_{\m}}=\a_{n}$
and
$\a'_{n_1+\cdots+n_{\m-1}+1},\ldots,\a'_{n_1+\cdots+n_{\m}}=\a'_{n}$;
\item[(2)] $p_{2i-1}(1)\in\FF_1$ and $p_{2i}(1)\in\FF_2$, for each
$i=1,\ldots,n$.
\end{itemize}

It is clear that, because of (2), the homeomorphism $\wti\b$,
associated to a $2n$-string braid $\b$ defining a special plat
presentation, keeps fixed both the sets $\FF_1$ and $\FF_2$.
Although a special plat presentation of a link is a very
particular case, we shall prove that every link admits such kind
of presentation.

\begin{proposition} \label{Proposition special} Every link $L$ admits a special
$2n$-plat presentation, for each $n\ge b(L)$.
\end{proposition}

\begin{proof}
Let $L$ be presented by a $2n$-plat. We show that this
presentation is equivalent to a special one, by using a finite
sequence of moves on the plat presentation which changes neither
the link type nor the number of plats. The moves are of the four
types $I$, $I'$, $II$ and $II'$ depicted in Figure 2. First of
all, it is straightforward that condition (1) can be satisfied by
applying a suitable sequence of moves of type $I$ and $I'$.
Furthermore, condition (2) is equivalent to the following: $(2')$
there exists an orientation of $L$ such that, for each
$i=1,\ldots,n$, the top arc $\a_i$ is oriented from $A_i$ to $B_i$
and the bottom arc $\a'_i$ is oriented from $B'_i$ to $A'_i$.
Therefore, choose any orientation on $L$ and apply moves of type
$II$ (resp. moves of type $II'$) to the top arcs (resp. bottom
arcs) which are oriented from $B_i$ to $A_i$ (resp. from $A'_i$ to
$B'_i$).
\end{proof}

%\bigskip

\begin{figure}[bht]
 \begin{center}
 \includegraphics*[totalheight=12cm]{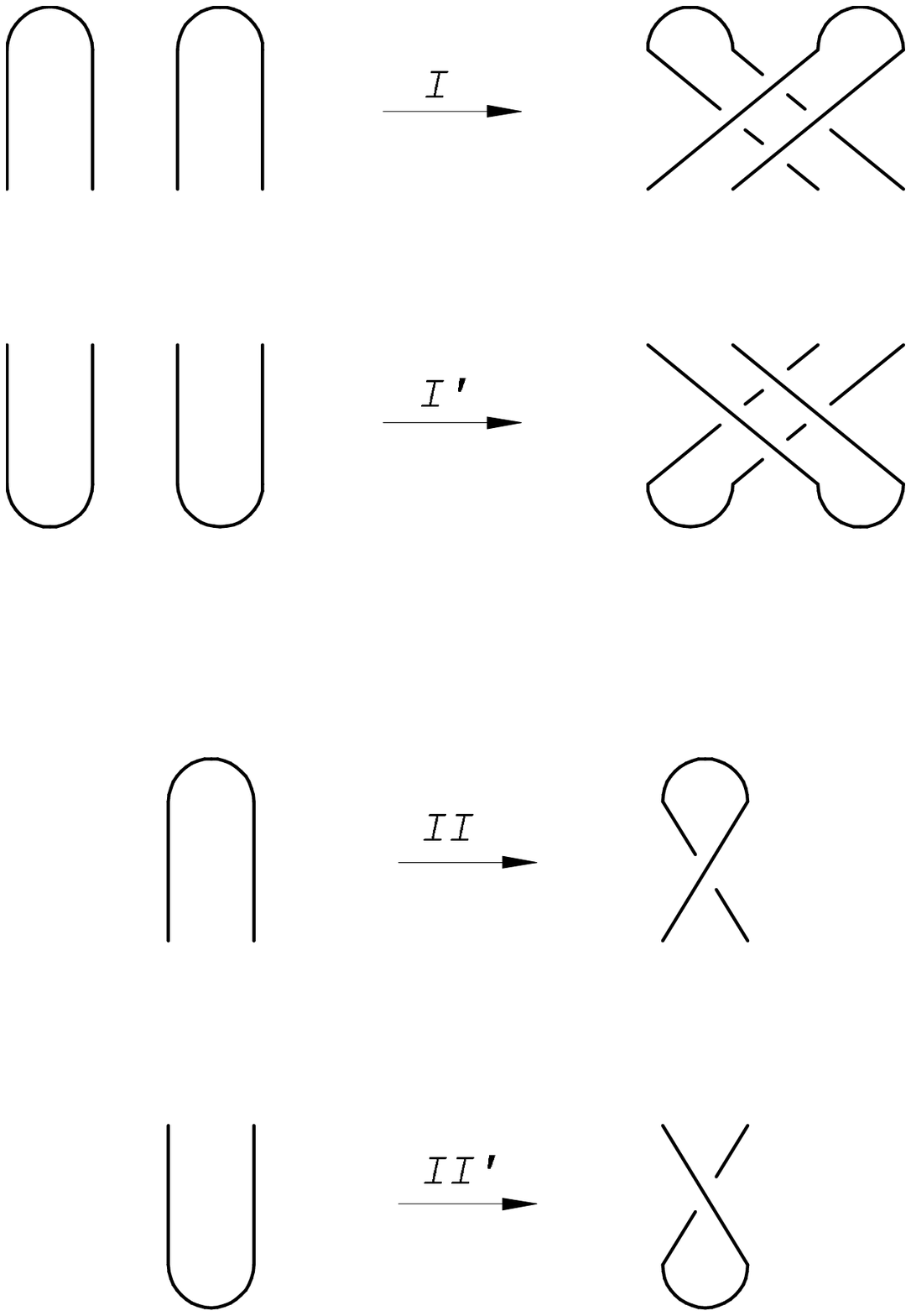}
 \end{center}
 \caption{Moves on plat presentations.}

 \label{Fig. 2}

\end{figure}

\medskip

A $p$-fold branched cyclic covering of an oriented $\m$-component
link $L=\bigcup_{j=1}^{\m}L_j\subset\S^3$ is completely determined
(up to equivalence) by assigning to each component $L_j$ an
integer $c_j\in{\bf Z}_p-\{0\}$, such that the set
$\{c_1,\ldots,c_{\m}\}$ generates the group ${\bf Z}_p$. The
monodromy associated to the covering sends each meridian of $L_j$,
coherently oriented with the chosen orientations of $L$ and
$\S^3$, to the permutation $(1\,2\,\cdots\,p)^{c_j}\in\Si_p$. By
multiplying each $c_j$ by the same invertible elements $c$ of
${\bf Z}_p$, we get an equivalent covering.

Following \cite{MM} we shall call a branched cyclic covering:
\begin{itemize}
\item [a)] {\it strictly-cyclic\/} if $c_{j'}=c_{j''}$, for every
$j',j''\in\{1,\ldots,\m\}$,
\item [b)] {\it almost-strictly-cyclic\/}
if $c_{j'}=\pm c_{j''}$, for every
$j',j''\in\{1,\ldots,\m\}$,\footnote{By a suitable reorientation
of the link, an almost-strictly-cyclic covering becomes a
strictly-cyclic one.}
\item [c)] {\it meridian-cyclic\/} if $\gcd(b,c_j)=1$, for every
$j\in\{1,\ldots,\m\}$,
\item [d)] {\it singly-cyclic\/} if $\gcd(b,c_j)=1$, for some
$j\in\{1,\ldots,\m\}$,
\item [e)] {\it monodromy-cyclic\/} if it is cyclic.
\end{itemize}

The following implications are straightforward: $$\text{ a)
}\Rightarrow\text{ b) }\Rightarrow\text{ c) }\Rightarrow\text{ d)
} \Rightarrow\text{ e) }.$$ Moreover, the five definitions are
equivalent when $L$ is a knot. Similar definitions and properties
also hold for a $p$-fold cyclic covering of a 3-ball, branched
over a set of properly embedded (oriented) arcs.

Note that every branched cyclic covering of a link arising from a
$p$-symmetric Heegaard splitting -- according to Birman-Hilden
construction -- is strictly-cyclic. We show that, conversely,
every $p$-fold branched strictly-cyclic covering of a link admits
a $p$-symmetric Heegaard splitting.

\begin{theorem} \label{Theorem 3} A $p$-fold strictly-cyclic
covering of $\S^3$ branched over a link $L$ of bridge number $b$
is a closed, orientable 3-manifold $M$ of $p$-symmetric Heegaard
genus $$g_p(M)\le(b-1)(p-1).$$
\end{theorem}

\begin{proof}
Let $L$ be presented by a special $2b$-plat arising from a braid
$\b$, and let $(\S^3,L)=(D,A_b)\cup_{\f'}(D',A'_b)$ be the
$(0,b)$-decomposition described in Remark 1. Now, all arguments of
the proofs of Theorem 3 of \cite{BH} entirely apply and the
condition of Lemma 4 of \cite{BH} is satisfied, since the
homeomorphism $\wti\b$ associated to $\b$ fixes both the sets
$\FF_1$ and $\FF_2$.
\end{proof}

\medskip

As a consequence of Theorem \ref{Theorem 3} and Birman-Hilden
results, there is a natural one-to-one correspondence between
$p$-symmetric Heegaard splittings and $p$-fold branched
strictly-cyclic coverings of links.

\section{Weakly $p$-symmetric Heegaard splittings}

A Heegaard splitting $M=Y_g\cup_{\ps}Y'_g$ of a 3-manifold $M$
will be called {\it weakly $p$-symmetric\/}, with $p>1$, if there
is an orientation-preserving homeomorphism $\R:Y_g\to Y_g$ of
period $p$ and  a homeomorphism $\r:Y_g\to Y'_g$ such that, if
$\HH$ is the cyclic group of order $p$ generated by $\R$ and
$\PS:\partial Y_g\to\partial Y_g$ is the orientation-preserving
homeomorphism $\PS=\r^{-1}_{\vert\partial Y'_g}\ps$, the following
conditions are fulfilled:
\begin{itemize}
\item [i')] $Y_g/\HH$ is homeomorphic to a 3-ball;
\item [ii')] $\bigcup_{h=1}^{p-1}\mbox{Fix}(\R^h)/\HH$ is an
unknotted set of arcs in the ball $Y_g/\HH$;
\item [iii')] there exists an integer $p_0$ such
that $\PS\R_{\vert\partial Y_g}\PS^{-1}=(\R_{\vert\partial
Y_g})^{p_0}$.
\end{itemize}

Observe that the map $\R'=\r\R\r^{-1}$ is obviously an
orientation-preserving homeomorphism of period $p$ of $Y'_g$ with
the same properties as $\R$, and the relation
$\ps\R_{\vert\partial Y_g}\ps^{-1}=(\R'_{\vert\partial
Y'_g})^{p_0}$ easily holds.

The {\it weakly $p$-symmetric Heegaard genus\/} ${\wti g}_p(M)$ of
a 3-manifold $M$ is the smallest integer $g$ such that $M$ admits
a weakly $p$-symmetric Heegaard splitting of genus $g$. Observe
that a $p$-symmetric Heegaard splitting, as defined in Section 1,
is weakly $p$-symmetric (take $\R=\P_{\vert Y_g}$ and $\r=\t$).

In order to obtain results analogous to the ones for $p$-symmetric
Heegaard splittings, for weakly $p$-symmetric Heegaard splittings,
we make use of some classical results on branched cyclic coverings
of $\S^2$. A $p$-fold cyclic covering $\g:M^2\to\S^2$, branched
over a set of $N$ points $B_{\g}=\{P_1,\ldots,P_N\}$ is completely
determined, up to equivalence, by a map
$\o'_{\g}:B_{\g}\to\Z_p-\{0\}$, $\o'_{\g}(P_k)=c_k$ such that (A)
$\{c_1,\ldots,c_N\}$ generates the group $\Z_p$ and (B)
$\sum_{k=1}^N c_k=0$. The monodromy map
$\o_{\g}:\p_1(\S^2-B_{\g})\to\Si_p$ associated to the covering
sends a loop $m_k$ represented by a small circle around $P_k$,
coherently oriented with a chosen orientation of $\S^2$, to the
permutation $(1\,2\,\cdots\,p)^{c_k}$, for all $k=1,\ldots,N$.
Condition (A) guarantees that $M^2$ is connected and condition (B)
depends on the fact that $\p_1(\S^2-B_{\g})$ admits the finite
presentation $$\p_1(\S^2-B_{\g})=<m_1,\ldots,m_N\mid \prod_{k=1}^N
m_k=1>.$$ The set $\g^{-1}(P_k)$ has cardinality $\gcd(p,c_k)$ for
all $k=1,\ldots,N$ and therefore, by standard calculations, $M^2$
has Euler characteristic:
$$\chi(M^2)=2p-Np+\sum_{k=1}^N\gcd(p,c_k).\eqno(**)$$

\medskip

\begin{theorem} \label{Theorem 2G} Every closed, orientable
3-manifold of weakly $p$-symmetric Heegaard genus $g$ admits a
representation as a $p$-fold cyclic covering of $\S^3$, branched
over a link $L$ of bridge number $$b(L)\le\frac{p-1+g}{p-p^*},$$
where $p^*=\max\{d\mid 1\le d<p\,,\,d\mbox{ divides } p\}$.
\end{theorem}
\begin{proof}
The proof is similar to the one of Theorem 2 of \cite{BH}. The
definition of weakly $p$-symmetric splitting implies that the
canonical projections $\p:Y_g\to Y_g/\HH$ and $\p':Y'_g\to
Y_g/\HH'$, where $\HH'$ is the cyclic group of order $p$ generated
by $\R'$, are $p$-fold branched cyclic coverings. By iii'), there
is a map $\ps':\partial Y_g/\HH\to\partial Y'_g/\HH'$ such that
$\p'_{\vert\partial Y'_g}\ps=\ps'\p_{\vert\partial Y_g}$.
Therefore, the map $\p\cup\p':Y_g\cup_{\ps}Y'_g\to
(Y_g/\HH)\cup_{\ps'}(Y'_g/\HH')\cong\S^3$ is a $p$-fold branched
cyclic covering. The restriction map $\g=\p_{\vert\partial Y_g}$
turns out to be a $p$-fold cyclic covering of $\partial
Y_g/\HH\cong\S^2$, branched over $2n$ points $P_1,\ldots,P_{2n}$,
with $\o'_{\g}(P_k)=c_k$ such that $c_{2h}=-c_{2h-1}$ for each
$h=1,\ldots,n$. By $(**)$ we get
$2-2g=2p-2np+2\sum_{h=1}^n\gcd(p,c_{2h-1})$. Since $c_k\ne 0$, we
have $(p,c_k)\le p^*$ and therefore $p-1+g=np-\sum_{h=1}^n
\gcd(p,c_{2h-1})\ge np-np^*$. Thus, $b(L)\le n\le(p-1+g)/(p-p^*)$
and the statement is achieved.
\end{proof}

\begin{corollary} \label{Corollary prime} If $p$ is prime,
then every closed, orientable 3-manifold of weakly $p$-symmetric
Heegaard genus $g$ admits a representation as a $p$-fold cyclic
covering of $\S^3$, branched over a link $L$ of bridge number
$$b(L)\le 1+\frac{g}{p-1}.$$
\end{corollary}
\begin{proof}
Straightforward, since $p^*=1$ when $p$ is prime.
\end{proof}

In order to prove the converse of Theorem \ref{Theorem 2G}, we
need the following result (compare Lemma 4 of \cite{BH}):

\begin{lemma} \label{Lemma 4W} Let $\g:M^2\to\S^2$ be a $p$-fold
cyclic covering, branched over the set
$B_{\g}=\{P_1,\ldots,P_N\}$. If $\PS':\S^2\to\S^2$ is an
orientation-preserving homeomorphism such that
$\PS'(B_{\g})=B_{\g}$ and $\o'_{\g}(\PS'(P_k))=\o'_{\g}(P_k)$, for
every $k=1,\ldots,N$, then $\PS'$ lifts to an
orientation-preserving homeomorphism $\PS:M^2\to M^2$ such that
$\PS'\g=\g\PS$.
\end{lemma}
\begin{proof}
It is well known that $\PS'$ lifts if and only if the induced
homomorphism $\PS'_*$ on the fundamental group $\p_1(\S^2-B_{\g})$
leaves the subgroup $H$ of the covering invariant \cite{Ma}. The
subgroup $H$ is the kernel of the homomorphism
$\o:\p_1(\S^2-B_{\g})=<m_1,\ldots,m_N\mid \prod_{k=1}^N
m_k=1>\to\Z_p$, defined by $\o(m_k)=c_k$, for each $k=1,\ldots,N$.
Then we have
$\o\PS'_*(m_k)=\o'_{\g}(\PS'(P_k))=\o'_{\g}(P_k)=\o(m_k)$ and
therefore $\o\PS'_*=\o$.
\end{proof}

%\medskip

\begin{theorem} \label{Theorem 3G} A $p$-fold cyclic covering
of $\S^3$, branched over a link $L$ of bridge number $b$, is a
closed, orientable 3-manifold $M$ of weakly $p$-symmetric Heegaard
genus $${\wti g}_p(M)\le(b-1)(p-1).$$
\end{theorem}
\begin{proof} Let $q:M\to\S^3$ be a $p$-fold cyclic
covering, branched over $L$, and let $L=\bigcup_{j=1}^{\m}L_j$ be
presented by a special $2b$-plat associated to a braid $\b$. If
$(\S^3,L)=(D,A_b)\cup_{\ps'}(D',A'_b)$ is the
$(0,b)$-decomposition of $L$ described in Remark 1, the map
$\r':D\to D'$ defined by $\r'(\infty)=\infty$ and
$\r'(x,t)=(x,1-t)$, for every $x\in\E^2$ and $t\le 0$, is an
homeomorphism which sends $\a_i$ onto $\a'_i$, for each
$i=1,\ldots,b$. Now, orient $L$ in such a way that, for each
$i=1,\ldots,b$, the arc $\a_i$ is oriented from $A_i$ to $B_i$ and
the arc $\a'_i$ is oriented from $B'_i$ to $A'_i$ (this is
possible since the plat is special). Let
$c_1,\ldots,c_{\m}\in\Z_p-\{0\}$ be the integers associated to the
components of $L$, according to the chosen orientations, and
defining the covering. For each $j=1,\ldots,\m$, the component
$L_j$ contains the $b_j$ top arcs
$\a_{b_1+\cdots+b_{j-1}+1},\ldots,\a_{b_1+\cdots+b_{j}}$ and the
$b_j$ bottom arcs
$\a'_{b_1+\cdots+b_{j-1}+1},\ldots,\a'_{b_1+\cdots+b_{j}}$. Let
$\p:Y_g\to D$ be the $p$-fold cyclic covering of the 3-ball $D$,
branched over the set of arcs $A_b$, with associated integers
$c_{j_1},\ldots,c_{j_b}$, where $L_{j_i}$ is the component of $L$
containing $\a_i$, for $i=1,\ldots,b$. If $Y'_g$ is another
handlebody of genus $g$ and $\r:Y_g\to Y'_g$ is a fixed
homeomorphism, then the map $\p'=\r'\p\r^{-1}:Y'_g\to D'$ is a
$p$-fold cyclic covering of $D'$, branched over $A'_b$, with
associated integer list $c_{j_1},\ldots,c_{j_b}$. Now let
$\g=\p_{\vert\partial Y_g}:\partial Y_g\to
\partial D$ be the restrictions of $\p$ to the surface $\partial Y_g$.
The map $\g$ is a $p$-fold cyclic covering of $\partial
D\cong\S^2$, branched over the $2b$ points $A_i,B_i\in\partial D$,
for $i=1,\ldots,b$, such that
$\o'_{\g}(A_i)=-\o'_{\g}(B_i)=c_{j_i}$. Because of properties (1)
and (2) of special plats, the homeomorphism
$\PS'=(\r'_{\vert\partial D})^{-1}\ps'$ of $\partial D$ satisfies
the condition of Lemma \ref{Lemma 4W}. Hence, $\PS'$ lifts to a
homeomorphism $\PS$ of $\partial Y_g$ such that $\PS'\g=\g\PS$. If
we define the homeomorphism $\ps=\r_{\vert\partial
Y_g}\PS:\partial Y_g\to\partial Y'_g$, the manifold
$Y_g\cup_{\ps}Y'_g$ turns out to be, by construction, a $p$-fold
branched cyclic covering of $L$, with the same monodromy as the
covering $q$ and therefore $Y_g\cup_{\ps}Y'_g$ is homeomorphic to
$M$. Moreover, $Y_g\cup_{\ps}Y'_g$ is a $p$-symmetric Heegaard
splitting of genus $g$. Indeed, let $\HH$ be the group of covering
transformations of $\p$ and $\R$ be a generator of $\HH$. Then
condition i') easily holds. As regards condition ii') we have:
$x\in\bigcup_{h=1}^{p-1}\mbox{Fix}(\R^h)/\HH\Leftrightarrow
\vert[x]_{\HH}\vert<p\Leftrightarrow\p(x)\in B_{\p}$. As for
condition iii'), the map $\PS\R_{\vert\partial Y_g}\PS^{-1}$ is
the lifting, with respect to $\g$, of
$\PS'(\PS')^{-1}=\mbox{Id}_{\partial D}$ and therefore it is a
covering transformation of $\p$. This proves point iii'). As far
as the genus of the splitting is concerned, we have by $(**)$:
$2-2g=2p-2bp+2\sum_{h=1}^b\gcd(p,c_{2h-1})$. Therefore ${\wti
g}_p(M)\le g=1-p+bp-\sum_{h=1}^b\gcd(p,c_{2h-1})\le
1-p+bp-b=(b-1)(p-1)$.
\end{proof}

\medskip

The case of branched cyclic coverings of two-bridge knots or links
is of particular interest.

\begin{corollary} \label{Corollary 2-bridge} A $p$-fold branched
cyclic covering of a two-bridge knot or link is a closed,
orientable 3-manifold $M$ of weakly $p$-symmetric Heegaard genus
$${\wti g}_p(M)\le p-1.$$
\end{corollary}
\begin{proof}
Follows from previous theorem, with $b=2$.
\end{proof}

\vspace{15 pt} {MICHELE MULAZZANI, Department
of Mathematics, University of Bologna, I-40127 Bologna, ITALY,
and C.I.R.A.M., Bologna, ITALY. E-mail: mulazza@dm.unibo.it}

\end{document}